\documentclass[a4]{article}
\font\est=eufb10 scaled\magstep1
\font\es=eufm10

\def\gg{\mbox{\es {g}}}
\def\gh{\mbox{\es {h}}}

\def\su{\mbox{\es{su}}}
\def\sp{\mbox{\es {sp}}}
\def\CC{\mbox{\es {C}}} 
\def\CCC{\mbox{\est {C}}}

\def\tr{\mbox{\rm {tr}}}
\def\det{\mbox{\rm {det}}}
\def\diag{\mbox{\rm {diag}}}
\def\Iso{\mbox{\rm {Iso}}}
\def\Hom{\mbox{\rm {Hom}}}
\def\R{\mbox{\boldmath $R$}}
\def\C{\mbox{\boldmath $C$}}
\def\H{\mbox{\boldmath $H$}}
\def\Z{\mbox{\boldmath $Z$}}      
\def\G{\mbox{\boldmath $G$}}

\def\0{\mbox{\boldmath {0}}}    
\def\1{\mbox{\boldmath {1}}}      
\def\2{\mbox{\boldmath {2}}}      
\def\3{\mbox{\boldmath {3}}}      
\def\4{\mbox{\boldmath {4}}}      
\def\5{\mbox{\boldmath {5}}}      
\def\6{\mbox{\boldmath {6}}}      
\def\7{\mbox{\boldmath {7}}}      
\def\8{\mbox{\boldmath {8}}}      
\def\9{\mbox{\boldmath {9}}}      
\def\a{\mbox{\boldmath $a$}}
\def\b{\mbox{\boldmath $b$}}

\def\e{\mbox{\boldmath $e$}}
\def\m{\mbox{\boldmath $m$}}

\def\m{\mbox{\boldmath $m$}}
\def\n{\mbox{\boldmath $n$}}

\def\dfrac#1#2{\displaystyle \frac{#1}{#2}}   
\begin{document}
\baselineskip=14pt
\title {\bf Adjoint orbit types of compact exceptional Lie group $G_2$ in its Lie algebra}   
\author         {Takashi M{\scriptsize IYASAKA}}
\date{}
\maketitle
\vspace{2mm}

{\large{\bf Introduction}} 

A Lie group $G$ naturally acts on its Lie algebra $\gg$, called the adjoint action. In this paper, we determine the orbit types of the compact exceptional Lie group $G_2$ in its Lie algebra $\gg_2$.  As results, the group $G_2$ has four orbit types in the Lie algebra $\gg_2$ as
$$
          G_2/G_2, \quad G_2/(U(1) \times U(1)), \quad G_2/((Sp(1)\times U(1))/\Z_2), \quad 
          G_2/((U(1)\times Sp(1))/\Z_2). $$
These orbits, especially the last two orbits, are not equivalent, that is, there exists no $G_2$-equivariant homeomorphism among them. Finally, the author would like to thank Professor Ichiro Yokota for his earnest guidance, useful advice and constant encouragement.		  
\vspace{4mm}

{\large{\bf 0. Preliminaries and notation}}
\vspace{3mm}
 
(1) \, For a group $G$ and  an element $s$ of $G$, $\widetilde{s}$ denotes the inner automorphism induced by $s:\widetilde{s}(g)=sgs^{-1}$, $g \in G$, then $G^{\widetilde{s}}= \{g \in G \, | \, sg = gs \}.$ Hereafter $G^{\widetilde{s}}$ is briefly written by $G^s$. 

(2) For a transformation group $G$ of a space $M$, the isotropy subgroup of $G$ at a point $m \in M$ is denoted by $G_m : G_m = \{g \in G \, | \, gm = m \}$.

(3) As mentioned in the introduction, a Lie group $G$ acts on its Lie algebra $\gg$. When $G$ is a compact Lie group, any element $X \in \gg$ is transformed to some element $D$ of a fixed Cartan subalgebra $\gh$. Hence, to determine the conjugate classes of isotropy subgroups $G_X$, it suffices to consider the case of $X = D \in \gh$. 
\vspace{4mm}

{\large{\bf 1. The Cayley alebra $\CCC$ and the group $G_2$}}
\vspace{3mm}

Let $\CC$ be the division Cayley algebra with the canonical basis $\{e_0=1,e_1,\dots ,e_7\}$ and in $\CC$ the conjugation $\overline{x}$, the inner product $(x, y)$ and the length $|x|$ are naturally defined ([1]).
The Cayley algebra $\CC$ contains naturally the field $\R$ of real numbers, furthermore the field $\C$ of complex numbers and the field $\H$ of quaternion numbers :
$$
\begin{array}{c}
    \R = \{x = x1 \, | \, x \in \R \}, \quad  \C = \{x_0 + x_1\e_1 \, | \, x_0, x_1 \in \R \}, \\
    \H = \{x_0 + x_1e_1 + x_2e_2 + x_3e_3 \, | \, x_0, x_1,x_2, x_3 \in \R \}.$$\end{array} $$

\noindent The automorphism group of the Cayley algebra $\CC$ :
$$
 G_2 = \{\alpha \in \Iso_R(\CC) \, | \, \alpha(xy) = (\alpha x)(\alpha y) \} $$
is the simply connected compact Lie group of type $G_2$.

\noindent Any element $x \in \CC$ can be expressed as $x = a + be_4, a, b \in \H,$ and $\CC$ is isomorphic to the algebra $\H \oplus \H e_4$ with multiplication
$$
        (a + be_4)(c + de_4) = (ac - \overline{d}b) + (b\overline{c} + da)e_4. $$     

\noindent We define an $\R$-linear transformation $\gamma$ of $\CC$ by
$$
      \gamma(a + be_4) = a - be_4, \quad a + be_4 \in \H \oplus \H e_4=\CC. $$ 
\noindent Next, to an element 
$$
           x = a + m_1e_2 + m_2e_4 + m_3e_6, \quad a, m_1, m_2, m_3 \in \C $$
       
\noindent of $\CC$, we associate an element
$$
                 a + \pmatrix{m_1 \cr
                              m_2 \cr
                              m_3} $$
of the algebra $\C \oplus \C^3$ with the multiplication
$$
        (a + \m)(b + \n) = (ab - \langle \m, \n \rangle)
                         + (a\n + \overline{b}\m - \overline{\m \times \n}), $$
\noindent where $\langle \m, \n \rangle = {}^t\m\overline{\n}$ and $\m \times \n \in \C^3$ is the exterior product of $\m, \n$. Note that $\C \oplus \C^3$ is a left $\C$-module. Hereafter we identify $\CC$ with $\H \oplus \H e_4$ and $\C \oplus \C^3$ :
$$
       \CC = \H \oplus \H e_4, \quad \CC = \C \oplus \C^3. $$

\noindent We define an $\R$-linear transformation $\gamma_1$ of $\CC$ by 
$$
      \gamma_1(a + \m) = \overline{a} + \overline{\m}, \quad a+\m \in  \C \oplus \C^3=\CC.$$
     
\noindent Then, $\gamma, \gamma_1 \in G_2, \gamma^2 = {\gamma_1}^2 = 1$ and $\gamma, \gamma_1$ are conjugate in $G_2$ ([1]).
\vspace{4mm}

{\large{\bf 2. Subgroups $(Sp(1) \times Sp(1))/\Z_2$ and $SU(3)$ of $G_2$}}
\vspace{3mm}

{\bf Proposition 1} ([1], [3]).\hspace{2mm}  $(G_2)^{\gamma} \cong (Sp(1) \times Sp(1))/\Z_2, \, \Z_2 = \{(1, 1), (-1, -1) \}.$
\vspace{2mm}

 {\bf Proof.} Let $Sp(1)=\{ p \in \H \, | \, p \overline{p}=1 \}$ and we define a map $\varphi : Sp(1) \times Sp(1) \to (G_2)^{\gamma}$ by 
 $$
      \varphi(p, q)(a + be_4) = qa\overline{q} + (pb\overline{q})e_4, \quad a + be_4 \in \H \oplus \H e_4 = \CC. $$
\noindent$\varphi$ is well-defined:  $\varphi(p, q) \in (G_2)^{\gamma}$ and $\varphi$ is a homomorphism. We shall show that $\varphi$ is onto. Note that
$$
      \CC_{\gamma} = \{x \in \CC \, | \, \gamma x = x \} = \H, \quad
      \CC_{- \gamma} = \{x \in \CC \, | \, \gamma x = - x \} = \H e_4 $$
\noindent and these spaces are invariant under the action of the group $(G_2)^{\gamma}$. Now, since $\alpha \in (G_2)^{\gamma}$ induces an atomorphism of $\CC_{\gamma} = \H$, there exists $q \in Sp(1)$ such that
$$
           \alpha a = qa\overline{q}, \quad a \in \H. $$
\noindent Let $\beta = \varphi(1, q)^{-1}\alpha$. Then $\beta \in (G_2)^{\gamma}$ and $\beta|\H = 1$. Since $\beta$ induces an endomorphism of $\H e_4$, there exists $p \in \H$ such that $\beta e_4 = pe_4$. From $|p| = |pe_4| = |\beta e_4| = |e_4| = 1$, we see that $p \in Sp(1)$. Then
$$
  \beta(a + be_4) = \beta a + (\beta b)(\beta e_4) = a + b(pe_4) 
                 = a + (pb)e_4 = \varphi(p, 1)(a + be_4). $$
\noindent Hence, $\beta = \varphi(p, 1)$ and $\alpha = \varphi(1, q)\varphi(p, 1) = \varphi(p, q)$ which shows that $\varphi$ is onto. Ker$\varphi = \{(1,1),(-1,-1)\}=\Z_2$. Thus, we have $(G_2)^{\gamma} \cong (Sp(1) \times Sp(1))/\Z_2$.
\vspace{3mm}

{\bf Proposition 2} ([3]).\hspace{1.5cm} $(G_2)_{e_1} \cong SU(3)$
\vspace{2mm}

{\bf Proof.} Let $SU(3)=\{ A \in M(3,\C)\, | \, A^*A=E,$ det$A=1 \}$ and we define a map $\psi : SU(3) \to (G_2)_{e_1}$, by 
$$
      \psi(A)(a + \m) = a + A\m, \quad a + \m \in \C \oplus \C^3=\CC. $$
$\psi$ is well-defined: $\psi(A) \in  (G_2)_{e_1}$. $\psi$ is injective and a homomorphism. We shall show that $\psi$ is onto. Let $\alpha \in (G_2)_{e_1}$. Note that $\alpha$ induces a \C-linear transformation of $\C^3$. Now let $$
              \alpha e_2 = \a_1, \quad \alpha e_4 = \a_2, \quad
              \alpha e_6 = \a_3 $$
\noindent and consider a matrix $A = (\a_1, \a_2, \a_3) \in M(3, \C)$.  From $(\alpha e_2)(\alpha e_4) = \alpha(e_2e_4) = - \alpha(e_6)$, we have $\a_1\a_2 = -\a_3$, namely, $- \langle \a_1, \a_2 \rangle - \overline{\a_1 \times \a_2} = - \a_3$, then
$$
       \langle \a_1, \a_2 \rangle = 0, \quad
        \a_3 = \overline{\a_1 \times \a_2}. $$
\noindent Similarly, we have $\langle \a_2, \a_3 \rangle = \langle \a_3, \a_1 \rangle = 0$. Moerover, from $(\alpha e_k)(\alpha e_k) = \alpha(e_ke_k) = \alpha(-1) = -1$, we have $\langle \a_k, \a_k \rangle = 1$, hence, $A \in U(3)=\{ A \in M(3,\C)\, | \, A^*A=E\}$. Finally, \det$A = (\a_3, \a_1 \times \a_2) = (\a_3, \overline{\a}_3) = \langle \a_3, \a_3 \rangle = 1$ (where $(\a, \b)$ is the usual inner product in $\C^3$ : $(\a, \b) = {}^t\a\b$). Hence, we have $A \in SU(3)$ and $\psi(A) = \alpha$ which shows that $\psi$ is onto. Thus, we have $SU(3) \cong (G_2)_{e_1}$.  
\vspace{4mm}

{\large{\bf 3. Adjoint orbit types of $G_2$ in the Lie algebra $\gg_2$}}
\vspace{3mm}

The Lie algebra $\gg_2$ of the group $G_2$ is given by
$$
           \gg_2 = \{D \in \Hom_R(\CC) \, | \, D(xy)
                 = (Dx)y + x(Dy) \}. $$
The adjoint action of the group $G_2$ on the Lie algebra $\gg_2$ is given by
$$
       \mu : G_2 \times \gg_2 \to \gg_2, \quad \mu(\alpha, D) = \alpha D\alpha^{-1}. $$
Now, we shall determine the adjoint orbit types of the group $G_2$ in $\gg_2$. Since the group $G_2$ contains the subgroup $SU(3)$ (Proposition 2), the Lie algeba $\gg_2$ also contains the subalgebra $\su(3) = \{D \in M(3, \C) \, | \, D^* + D = 0, \tr(D) = 0 \}$. We choose a Cartan subalgebra $\gh$ of $\su(3)$ as
$$
     \gh = \bigg\{D(r, s, t) = \pmatrix{re_1 & 0 & 0 \cr
                                           0 & se_1 & 0 \cr                                                                0 & 0 & te_1} \left| \, 
            r, s, t \in \R, r + s + t = 0 \right. \bigg\}, $$
which is also a Cartan subalgebra of $\gg_2$. The order of $r, s, t$ of $D(r, s, t)$ is arbitrarily exchanged by the action of $G_2$. \\
Note that between the induced mappings $\varphi_* : \sp(1) \oplus \sp(1) \to \gg_2$ and $\psi_* :\su(3) \to \gg_2$ of $\varphi$ and $\psi$, there exist the following relations.
$$
      \varphi_*(e_1, 0) = \psi_*(\diag(0, e_1, -e_1)), \quad
      \varphi_*(0, e_1) = \psi_*(\diag(2e_1, -e_1, -e_1)).$$       
\vspace{3mm}

{\bf Theorem 3.} \, {\it The orbit types in $\gg_2$ through $D(r, s, t) \in \gh$ under the adjoint acion of the group $G_2$ are as follows.}

(1){\it In the case of $r = s = t = 0$, the orbit through $D(0, 0, 0)$ is}
$$
      G_2/G_2. $$ 

(2) {\it In the case of $r, s, t$ are non-zero and distinct, the orbit through $D(r, s, t)$ is}
$$
      G_2/(U(1) \times U(1)).  $$

(3) {\it In the case of $r$ is non-zero, the orbit through $D(2r, -r, -r)$ is}
$$
      G_2/((Sp(1) \times U(1))/\Z_2). $$

(4) {\it In the case of $r$ is non-zero, the orbit through $D(0, r, - r)$ is}
$$
      G_2/((U(1) \times Sp(1))/\Z_2). $$
	  
{\bf Proof.} (1) is trivial.

(2) \, Let $U(1)= \{a \in \C \, | \, a\overline{a}=1 \}$ and $S(U(1) \times U(1) \times U(1))$ be the diagonal subgroup of $SU(3)$ and $\varphi : S(U(1) \times U(1) \times U(1)) \to (G_2)_{D(r, s, t)}$ be the restriction map $\varphi$ of Proposition 1. Then, $\varphi$ is well-defined and injective. We shall show that $\varphi$ is onto. For this purpose, we first show that for $\alpha \in (G_2)_{D(r, s, t)}$, we have
$$
       \alpha e_1 = \pm e_1.$$
Indeed, let $\alpha e_1 = a + \m \in \C \oplus \C^3$. From the condition $\alpha D(r, s, t) = D(r, s, t)\alpha$,
$$
    0 = \alpha 0 = \alpha D(r, s, t)e_1 = D(r, s, t)\alpha e_1 
      = D(r, s, t)(a + \m) = D(r, s, t)\m, $$
we have $\m = 0$. So $\alpha e_1 = a$. Since $|a| = |\alpha e_1| = |e_1| = 1$ and $\alpha 1 = 1$, we have $\alpha e_1 = \pm e_1$.    

\noindent In the case of $\alpha e_1 = - e_1$, consider $\gamma_1 \in G_2$. Then $\alpha\gamma_1e_1 = e_1$, so $\alpha\gamma_1 = A = \Big( a_{kl} \Big) \in SU(3)$ (Proposition 2). Hence, $\alpha = A\gamma_1$. Then   
$$
     A\gamma_1D(r, s, t)\pmatrix{1 \cr 0 \cr 0}
         = A\gamma_1\pmatrix{re_1 \cr 0 \cr 0}
         = A\pmatrix{- re_1 \cr 0 \cr 0}
         = \pmatrix{- ra_{11}e_1 \cr - ra_{21}e_1 \cr - ra_{31}e_1}.   $$
On the other hand,
$$
     D(r, s, t)A\gamma_1\pmatrix{1 \cr 0 \cr 0}
         = D(r, s, t)A\pmatrix{1 \cr 0 \cr 0}
         = D(r, s, t)\pmatrix{a_{11} \cr a_{21} \cr a_{31}}
         = \pmatrix{re_1a_{11} \cr se_1a_{21} \cr te_1a_{31}}. $$
From the condition $\alpha D(r, s, t) = D(r, s, t)\alpha$, we have
$$
   {}^t(- ra_{11}e_1, - ra_{21}e_1,\\ - ra_{31}e_1) = {}^t(re_1a_{11}, se_1a_{21}, te_1a_{31}). $$
 But, it is easy to see that this is false. Hence, we have $\alpha e_1 = e_1$, so $\alpha = A \in SU(3)$. From the condition $\alpha D(r, s, t) = D(r, s, t)\alpha$ again, we have $\alpha = A \in S(U(1) \times U(1) \times U(1))$, which shows that $\varphi$ is onto. Thus we have $(G_2)_{D(r, s, t)} = S(U(1) \times U(1) \times U(1)) \cong U(1) \times U(1)$.
\vspace{1mm}

(3) \, Since $\exp(\pi D(2, -1, -1)) = \exp(\pi\psi_*(\diag(2e_1, -e_1, -e_1)))= \exp(\pi\varphi_*(0,$ $ -e_1)) = \gamma,$ if $\alpha \in G_2$ commutes with $D(2, -1, -1)$, then $\alpha$ also commutes with $\gamma$. Hence, $ \alpha \in (G_2)^{\gamma} = \varphi(Sp(1) \times Sp(1))$ (Proposition 1). So there exist $p, q \in Sp(1)$ such that $\alpha =$ 
\vspace{1mm}

\noindent $\varphi(p, q)$. Again from the commutativity with $\exp(\dfrac{\pi}{2}D(2, -1, -1)) = \exp(\dfrac{\pi}{2}\psi_*(\diag(0, 0, e_1)))$\vspace{1mm}

\noindent $= \varphi(1, e_1)$, we have $\varphi(p, q)\varphi(1, e_1)$ $ = \varphi(1, e_1)\varphi(p, q)$, hence, $\varphi(p, qe_1) = \varphi(p, e_1q)$. So $qe_1 = e_1q$, therefore  $q\in \C \cap Sp(1) = U(1)$. Conversely, $\alpha = \varphi(p, q) (p \in Sp(1),q \in U(1))$ commutes with $\varphi(1, e^{te_1})$, $t \in \R$, so $\alpha$ also commutes with $\varphi_{*}(0, e_1) = D(2, -1, -1)$. Thus, we have $(G_2)_{D(2, -1, -1)} = \varphi(Sp(1) \times U(1)) = (Sp(1) \times U(1))/\Z_2.$      
\vspace{1mm}                                   

(4) $\exp(\pi D(1, -1, 0)) = \exp(\pi\psi_*(\diag(e_1, -e_1, 0))) = \exp(\pi(\varphi_*(e_1, 0)) = \gamma$. Hence, $(G_2)_{D(r, -r, 0)} = (U(1) \times Sp(1))/\Z_2$ is proved in a similar way to (2) above.

\vspace{4mm}

{\bf Proposition 4.} {\it The groups $(G_2)_{D(0,r, -r)} \cong(U(1) \times Sp(1))/\Z_2$ and 
\\ $(G_2)_{D(2s, -s, -s)} \cong (Sp(1) \times U(1))/\Z_2$ are not conjugate in the group $G_2$.}
\vspace{2mm}

{\bf Proof.} We shall prove that $D(0, r, -r)$ and $D(2s, -s, -s)$ are not conjugate under the action of the group $G_2$. Suppose that there exists $\alpha \in G_2$ such that
$$
        \alpha (D(0, r, -r)) = (D(2s, -s, -s))\alpha. $$
Let $\alpha e_2=a_0 + \m \in \C \oplus \C^3$, $\m ={}^t(m_1, m_2, m_3)$. Then, from
$$
        \alpha\pmatrix{0 & 0 & 0\cr
                      0 & re_1 & 0 \cr
                      0 & 0 & -re_1}\pmatrix{1 \cr 0 \cr 0}
        = \pmatrix{2s e_1& 0 & 0 \cr
                  0 & -se_1 & 0 \cr
                  0 & 0 & -se_1}\Big(a_0 + \pmatrix{m_1 \cr m_2 \cr m_3}\Big), $$
we have $\pmatrix{0 \cr 0 \cr 0} = \pmatrix{2sm_1e_1 \cr -sm_2e_1 \cr -sm_3e_1}$, hence,  $m_1 = m_2 = m_3 = 0$, so
$$
                \alpha e_2 = a_0 \in \C. $$
From $(\alpha e_2)(\alpha e_2) = \alpha(e_2e_2) = \alpha(-1) = -1$, we have $a_0a_0 = -1$, hence, $\alpha e_2 = a_0 = \pm e_1$. Similarly, we have $\alpha e_3 = \pm e_1$. Then $\alpha e_1 = (\alpha e_2)(\alpha e_3) = (\pm e_1)(\pm e_1) = \pm 1$, which is a contradiction.

\vspace{4mm}

{\bf Theorem 5.} {\it  The orbit spaces 
$X=\{\alpha (D(0,r,-r))\alpha^{-1}\, | \,\alpha \in G_2\}$ and $Y=\{\alpha (D(2s,\\
-s,-s))\alpha^{-1}\, | \, \alpha \in G_2 \}$ are not equivalent.}
\vspace{2mm}

{\bf Proof.} Suppose that there exists a $G_{2}$-equivariant homeomorphism $h:X \to Y$. Then, there exists $\alpha \in G_2$ such that 
\begin{equation}
h(D(0,r,-r))=\alpha(D(2s,-s,-s) )\alpha^{-1} .
\end{equation}
For any  $\delta \in (G_2)_{D(0,r,-r)}$, we have $\delta (h(D(0,r,-r))) \delta^{-1}=\delta (\alpha (D(2s,-s,-s)) \alpha^{-1})\delta^{-1}$.
Since h is a $\G_2-$equivariant homeomorphism, we have $\widetilde{\delta} h=h \widetilde{\delta}$. Hence, we see
\begin{equation} 
 h(D(0,r,-r))=\delta (\alpha (D(2s,-s,-s))\alpha^{-1})\delta^{-1}. 
 \end{equation}
 From (2) and (1), we have 
 $\delta (\alpha (D(2s,-s,-s) )\alpha^{-1})\delta^{-1} = \alpha (D(2s,-s,-s)) \alpha^{-1}$, that is, $\alpha^{-1}(\delta (\alpha (D(2s,-s,-s)) \alpha^{-1})\delta^{-1}) \alpha = D(2s,-s,-s) $. This implies $\alpha^{-1}\delta \alpha \in (G_2)_{D(2s,-s,-s)}.$ Thus, we have 
$$
\alpha^{-1}((G_2)_{D(0,r,-r)})\alpha \subset (G_2)_{D(2s,-s,-s)}.$$
Since dim$((G_2)_{D(0,r,-r)})=$dim$((G_2)_{D(2s,-s,-s)})$ (Theorem 3.(3) and (4)), the inclusion above must be equal, that is, 
$$
\alpha^{-1}((G_2)_{D(0,r,-r)})\alpha = (G_2)_{D(2s,-s,-s)}.$$
This means that the groups $(G_2)_{D(0,r,-r)}$ and $(G_2)_{D(2s,-s,-s)}$ are conjugate in $G_2$, which contradicts Proposition 4.

\bigskip
\begin{flushright}
\begin{tabular}{l}
{\sc Takashi Miyasaka} \\
{\sc Takatoh High School} \\
{\sc Obara, Takatoh, }\\
{\sc 396-0293, Japan }\\
E-mail: coolkai@mac.com
\end{tabular}
\end{flushright}

 \end{document}